\documentclass[12pt]{article}
\usepackage{amsmath}
\usepackage{amssymb}
\usepackage{graphicx}
\usepackage{xcolor}
\begin{document}

\title{Reduction of Stratified Axi-Symmetric Euler-Poisson 
Equations Under Symmetry} 
                   
\author{Mayer Humi\thanks {e-mail: mhumi@wpi.edu.}\\
Department of Mathematical Sciences\\
Worcester Polytechnic Institute\\
100 Institute Road\\
Worcester, MA  01609}

\maketitle
\thispagestyle{empty}

\begin{abstract}

The paper considers Euler-Poisson equations which govern the 
steady state of a self gravitating, rotating, axi-symmetric fluid under the 
additional assumption that it is incompressible and stratified. 
In this setting we show that the original system of six nonlinear partial 
differential equations can be reduced to 
two equations, one for the mass density and the other for gravitational field.
This reduction is carried out in cylindrical 
coordinates. As a result we are able to derive also expressions for the 
pressure as a function of the density. The resulting equations are then 
solved analytically. 
These analytic solutions are used then to determine the shape of the 
rotating star (or interstellar cloud) by applying the boundary condition 
that the pressure is zero at the boundary.
\end{abstract}

\newpage

\section{Introduction}

The steady states of self gravitating fluid in three dimensions 
have been studied by a long list of theoretical physicists and astrophysicists. 
(For an extensive list of references see \cite{SC1,SC2,MHT,EAM,JT,CSY}. 
In fact research related to this problem persists even today 
\cite{JJM,YYL,LS,AP1,AP2,IWR,VPS}. The motivation for 
this research is due to the interest in the formation, shape, and stability of 
stars and other celestial bodies.  

Within the context of classical mechanics attempts to describe star interiors
are based on Euler-Poisson equations \cite{SC1,SC2}. Well known solutions 
to these equations are the Lane-Emden functions, which describe static 
steady state of non-rotating spherically symmetric fluid with mass-density 
$\rho=\rho(r)$ and flow field $\bf u=0$.
The generalization of these equations to include axi-symmetric rotations was 
considered by Milne\cite{EAM}, Chandrasekhar\cite{SC1,SC2}, and 
many others \cite{KN,AK,EAM,OVM}. 
One of  difficulties in the treatment of this problem is due to the fact that 
the boundary of the domain can not be prescribed apriori, and one has to address a free boundary problem. An approximate treatment of this problem for 
polytropic stars in spherical coordinates was made in the seminal paper by
Roxburgh\cite{IWR}. Other treatments which considered different aspects
of this problem appeared in the literature since then 
\cite{VSB,AT,JJM,KNM,KN,AP2}.  

{In a previous paper \cite{MH3} on this topic, we addressed 
the steady states of non-rotating self gravitating incompressible fluid 
with axial symmetry. In the present paper, we generalize this treatment 
and address the modeling of axi-symmetric {\bf rotating} fluids}. To do so 
we add the assumption that the mass-density is stratified 
\cite{MH1,MH2,MH3,MH4,MH6,MHT,JT,CSY}
to the Euler-Poisson equations with axi-symmetric rotations. Under these 
assumptions we show that the number of model equations for the (non static)
steady state can be reduced from six to a system of two coupled equations. 
One for the mass-density and the second for the gravitational field. 
These equations 
contain, however, a parameter  function $h(\rho)$ that encode the information 
about the momentum distribution within the star. This reduction in the number 
of model equations (in this settings) may be used to obtain new 
insights for the treatment of this problem and make it tractable both 
analytically and numerically. We provide in this paper analytic and numerical
solutions to these equations and apply these solutions to solve for the shape 
of an axi-symmetric rotating star.

It might be argued that Euler-Poisson equations do not actually hold
in a star interior due to the various physical processes taking place
there (e.g turbulence, radiation, compressibility, etc)\cite{MH5}. 
Nevertheless they provide a natural extension to the results on the 
equilibrium states of three dimensional bodies under gravity.
  
The plan of the paper is as follows: In Sec. $2$, we present the basic 
model equations. In Sec. $3$, we carry out, in cylindrical coordinates,
the reduction of these model equations from six to two. (similar reduction 
can be carried out in spherical coordinates).
We provide also expressions for the pressure in this coordinate system.
In Secs. $4$ and $5$, we derive analytic solutions to these equations and 
discuss the shape of the rotating star (or interstellar cloud) by applying 
the boundary condition that the pressure is zero on the boundary. We end up 
in Sec. $6$, with summary and conclusions.

\setcounter{equation}{0}
\section{Derivation of the Model Equations}

In this paper we consider the state of an inviscid 
incompressible stratified self gravitating fluid. In addition we assume that 
the fluid it is subject to axial rotations. 
The hydrodynamic equations that govern this flow in an inertial frame of
reference are \cite{VSB,SC1,MH1,MH4,OVM,JT};
\begin{equation}
\label{2.1}
\nabla \cdot {\bf v} = 0
\end{equation}
\begin{equation}
\label{2.2}
{\bf v} \cdot \nabla\rho = 0
\end{equation}
\begin{equation}
\label{2.3}
\frac{1}{2}\rho\nabla({\bf v}\cdot{\bf v})+
\rho(\nabla\times{\bf v})\times {\bf v} = 
-\nabla p -\rho \nabla \Phi 
\end{equation}
\begin{equation}
\label{2.4}
\nabla^2 \Phi = 4 \pi G \rho
\end{equation}

where ${\bf v}=(u,v,w)$ is the fluid velocity, $\rho$ is its density 
$p$ is the pressure, $\Phi$ is the gravitational potential, G is 
the gravitational constant, and the momentum equations (\ref{2.3}) are 
written in Lambs's form. Subscripts denote differentiation with respect
to the indicated variable.  

We can nondimensionalize these equations by introducing the following scalings
\begin{eqnarray}
\label{2.5}
&&x= L\tilde{x},\,\,\ y=L\tilde{y},\,\,\, z=L\tilde{z},\,\,\, 
{\bf v}=U_0 \tilde{{\bf v}} ,\,\,\\ \notag
&&\rho = \rho_0 \tilde{\rho},\,\,\ p=\rho_0 U_0^2\tilde{p},\,\,\,
\Phi= U_0^2 \tilde{\Phi},\,\,\, \omega=\frac{U_0}{L}\tilde{\omega}.
\end{eqnarray}
where $L,U_0,\rho_0$ are some characteristic length, velocity, and mass density
respectively that characterize the problem at hand.

Substituting these scalings in (\ref{2.1})-(\ref{2.4}) and dropping 
the tildes, these equations remain unchanged (but the quantities that appear
in these equations become nondimensional) while $G$ is replaced by 
$\tilde{G}=\frac{G\rho_0 L^2}{U_0^2}$. (Once again we drop the tilde).

We now restrict our discussion to bodies which are axi-symmetric.
Without loss of generality, we shall assume henceforth that this axis
of symmetry coincides with the z-axis. Under this assumption, it is
expeditious to treat the flow in cylindrical (or spherical) coordinate
system. In standard cylindrical coordinates $(r,\theta,z)$, 
we then have (due the symmetry) ${\bf v}= {\bf v}(r,z)$,  i.e. the flow 
and the other functions that appear in (\ref{2.1})-(\ref{2.4}) are
independent of the angle $\theta$. 

{\bf We note that all functions in this paper will be assumed to be smooth. 
Furthermore for the rest of this paper we exclude from our discussion 
the exceptional cases where $\rho$ is constant or ${\bf v}={\bf 0}$. 
To determine the shape of a star under these assumptions we impose the 
boundary condition $p=0$ (\cite{JT} p. $54$ and p. $121$). In addition 
$\rho$ must (obviously) satisfy $\rho \ge 0$.}

\setcounter{equation}{0}
\section{Reduction in Cylindrical Coordinates}

Following the standard notation we introduce the frame
$$
{\bf e}_r =(\cos\theta,\sin\theta,0),\,\,\, {\bf e}_{\theta} =
(-sin\theta,\cos\theta,0),\,\,\, {\bf e}_z=(0,0,1).
$$
In this frame we have under present assumptions
\begin{equation}
\label{2.00}
{\bf v} = u(r,z){\bf e}_r + w(r,z){\bf e}_z+
v(r,z){\bf e}_{\theta}={\bf u}(r,z)+v(r,z){\bf e}_{\theta}
\end{equation}
The momentum equations for ${\bf u}$ can be written as
\begin{equation}
\label{2.02}
\rho{\bf u}\cdot \nabla {\bf u}=-\nabla p-
\rho\nabla\Phi+ \rho\frac{v^2}{r}{\bf e}_r.
\end{equation}
The equation for $v$ is
\begin{equation}
\label{2.03}
{\bf u}\cdot\nabla v+\frac{uv}{r}=0.
\end{equation}
We observe also that we can replace ${\bf v}$  by ${\bf u}$ in 
(\ref{2.1})-(\ref{2.2}).

In the cylindrical coordinate system the continuity equation (\ref{2.1}) 
becomes
\begin{equation}
\label{2.6}
\frac{1}{r}\frac{\partial(ru)}{\partial r} +\frac{\partial w}{\partial z}=0.
\end{equation}
This can be rewritten as
\begin{equation}
\label{2.7}
\frac{1}{r}\left[\frac{\partial(ru)}{\partial r} +\frac{\partial (rw)}{\partial z}\right]=0.
\end{equation}
It follows then that it is appropriate to introduce Stokes stream function 
$\psi$ \cite{IWR,JT,CSY}, which satisfy
\begin{equation}
\label{2.8}
u=\frac{1}{r}\frac{\partial \psi}{\partial z},\,\,\, w=-\frac{1}{r}\frac{\partial \psi}{\partial r}.
\end{equation}
Observe that since we excluded the case where ${\bf v}= {\bf 0}$, $\psi$
can not  be a constant. 
With these definitions (\ref{2.1}) is satisfied automatically by $\psi$. 

Since $\rho=\rho(r,z)$, (\ref{2.2}) in this frame is
\begin{equation}
\label{2.9}
u\rho_r+w\rho_z=0.
\end{equation}
Expressing $u,w$ in terms of $\psi$ we obtain
\begin{equation}
\label{2.10}
J\{\rho,\psi \}=0,
\end{equation}
where for any two (smooth) functions $F,G$
\begin{equation}
\label{2.11}
J\{F,G\}=\frac{\partial F}{\partial r}\frac{\partial G}{\partial z} -
        \frac{\partial F}{\partial z}\frac{\partial G}{\partial r}.
\end{equation}
{We observe that unless either $\rho$ or $\psi$
are constants, (\ref{2.10}) implies that we can express $\rho=\rho(\psi)$ or
$\psi=\psi(\rho)$. 

The explicit form of (\ref{2.03}) is
\begin{equation}
\label{2.05}
u\left(\frac{\partial v}{\partial r}+\frac{v}{r}\right)+
w\frac{\partial v}{\partial z}=0.
\end{equation}
Substituting $v =\frac{{\tilde v}}{r}$ this equation becomes
\begin{equation}
\label{2.04}
u{\tilde v}_r + w{\tilde v }_z=0
\end{equation}
Since \eqref{2.9} for $\rho$ and \eqref{2.04} for ${\tilde v}$ are the same 
it is natural to assume that there exists a smooth function $f$ so 
that ${\tilde v} =f(\rho)$. Hence
\begin{equation}
\label{2.04a}
v=\frac{f(\rho)}{r}
\end{equation}

{We observe that the function $f(\rho)$ should be 
determined by the physical constraints that are being imposed on the body 
by its rotation. "In absence" of such (explicit) constraints, it can be 
considered as a "parameter function"}

The momentum equations (\ref{2.02}) in this coordinate system become
\begin{equation}
\label{2.12}
\rho (uu_r +wu_z) = -p_r -\rho\Phi_r +\rho \frac{f(\rho)^2}{r^3}
\end{equation}
\begin{equation}
\label{2.13}
\rho (uw_r+ww_z) = -p_z -\rho\Phi_z,
\end{equation}

To eliminate $p$ from (\ref{2.12}) and  (\ref{2.13}) we differentiate 
these equations with respect to $z,r$ respectively and subtract. We obtain;
\begin{eqnarray}
\label{2.14}
&&\rho_r(uw_r+ww_z)+\rho(uw_r+ww_z)_r- \\ \notag
&&\rho_z(uu_r+wu_z)-\rho(uu_r+wu_z)_z= \\ \notag
&&-J\{\rho,\Phi\}-J\{\rho,\frac{H(\rho)}{r^2}\},
\end{eqnarray}
where
\begin{equation}
\label{2.14a}
H(\rho)=\frac{f^2}{2}+\rho f\,f_{\rho}.
\end{equation}
For the first and third terms on the left hand side of this equation, 
we obtain using (\ref{2.9})
\begin{eqnarray}
\label{2.15}
&&\rho_r(uw_r+ww_z)-\rho_z(uu_r+wu_z) = \\ \notag
&&\rho_r(ww_z+uu_z)-\rho_z(uu_r+ww_r)= \\ \notag
&&J\{\rho,\frac{u^2+w^2}{2}\}
\end{eqnarray}
Similarly for the second and forth terms on the left hand side of (\ref{2.14}),
we have
\begin{equation}
\label{2.16}
\rho\left[(uw_r+ww_z)_r-(uu_r+wu_z)_z\right] = \rho\left[u\chi_r + w\chi_z+(u_r+w_z)\chi \right]
\end{equation}
where $\chi=w_r-u_z$. However, from (\ref{2.6}) we have
$$
u_r+w_z= -\frac{u}{r}.
$$
Using this equality and expressing $u,w$ in terms of $\psi$ leads to
\begin{equation}
\label{2.17}
u\chi_r + w\chi_z+(u_r+w_z)\chi =-J\{\psi,\frac{\chi}{r}\}.
\end{equation}

Hence, we finally obtain that
\begin{equation}
\label{2.18}
\rho\left[(uw_r+ww_z)_r-(uu_r+wu_z)_z\right] = \rho J\left\{\psi,\frac{1}{r^2}\left(\nabla^2\psi-\frac{2}{r}\frac{\partial \psi}{\partial r}\right)\right\} .
\end{equation}
Combining the results of (\ref{2.14}),(\ref{2.15}) and (\ref{2.18}) it 
follows that
\begin{eqnarray}
\label{2.19}
&&J\{\rho,\frac{u^2+w^2}{2}\}+ \rho J\left\{\psi,\frac{1}{r^2}\left(\nabla^2\psi-\frac{2}{r}\frac{\partial \psi}{\partial r}\right)\right\}= \\ \notag
&&-J\{\rho,\Phi\}-J\{\rho,\frac{H(\rho)}{r^2}\}.
\end{eqnarray}
To express (\ref{2.19}) in terms of $\rho$ only, we use the fact that 
$\psi=\psi(\rho)$ and, therefore,
\begin{equation}
\label{2.20}
\psi_r =\psi_{\rho}\rho_r,\,\,\ \psi_z =\psi_{\rho}\rho_z,\,\,\
\nabla^2\psi=\psi_{\rho\rho}[\rho_r^2+\rho_z^2] + \psi_{\rho}\nabla^2\rho.
\end{equation}
Using these relations we have 
\begin{equation}
\label{2.21}
J\left\{\rho,\frac{u^2+w^2}{2}\right\} = J\left\{\rho,\frac{\psi_{\rho}^2}{2r^2}(\rho_r^2+\rho_z^2)\right\}
\end{equation}
\begin{equation}
\label{2.22}
\rho J\left\{\psi,\frac{1}{r^2}\left(\nabla^2\psi-\frac{2}{r}\frac{\partial\psi}{\partial r}\right)\right\} = \rho J\left\{\rho,\frac{1}{r^2}\left[\psi_{\rho}^2(\nabla^2\rho -\frac{2}{r}\rho_r) +\psi_{\rho}\psi_{\rho\rho}(\rho_r^2+\rho_z^2)\right]\right\}
\end{equation}
Substituting these results in (\ref{2.19}) leads to
\begin{equation}
\label{2.23}
J\left\{\rho,\frac{\rho}{r^2}\left[\psi_{\rho}^2(\nabla^2\rho -\frac{2}{r}\rho_r) +\psi_{\rho}\psi_{\rho\rho}(\rho_r^2+\rho_z^2)\right] +\frac{\psi_{\rho}^2}{2r^2}(\rho_r^2+\rho_z^2) +\Phi +\frac{H(\rho)}{r^2}\right\} = 0.
\end{equation}
{This is satisfied if there exists some function $S(\rho)$ 
such that}
\begin{equation}
\label{2.24}
\frac{\rho}{r^2}\left[\psi_{\rho}^2(\nabla^2\rho -\frac{2}{r}\rho_r) +
\psi_{\rho}\psi_{\rho\rho}(\rho_r^2+\rho_z^2)\right] +
\frac{\psi_{\rho}^2}{2r^2}(\rho_r^2+\rho_z^2) +\Phi+\frac{H(\rho)}{r^2} =S(\rho)
\end{equation}
{where $S(\rho)$ is some function of $\rho$ which can be viewed 
as a ``gauge". In the following we let $S(\rho)=0$.}

Introducing,
$$
h(\rho)=\rho\psi_{\rho}^2,\,\,\, h^{\prime}(\rho)=\frac{dh(\rho)}{d\rho}.
$$
We can rewrite (\ref{2.24}) more succinctly
\begin{equation}
\label{2.25}
h(\rho)\left(\nabla^2\rho -\frac{2}{r}\rho_r \right) +
\frac{h^{\prime}(\rho)}{2}(\rho_r^2+\rho_z^2) = r^2\left( S(\rho) -
\Phi-\frac{H(\rho)}{r^2}\right)
\end{equation}
This can be rewritten in the form
\begin{equation}
\label{2.29}
h(\rho)^{1/2}\nabla {\bf{\cdot}} ( h(\rho)^{1/2}\nabla \rho)-
\frac{2h(\rho)}{r}\rho_r = r^2\left(S(\rho) -\Phi-\frac{H(\rho)}{r^2}\right) .
\end{equation}
{\bf Thus, we reduced the original nonlinear system of six
partial differential equations (\ref{2.1})-(\ref{2.4}) to a coupled system
of two second order equations consisting of (\ref{2.4}) and (\ref{2.29})}.

One can use a transformation to simplify (\ref{2.29}) as
follows: First introduce $q(\rho)$ so that
\begin{equation}
\label{2.29a} 
\frac{dq}{d\rho}=h(\rho)^{1/2}.
\end{equation}
Hence, 
\begin{equation}
\label{2.29b}
\frac{\partial q}{\partial r}=h(\rho)^{1/2}\frac{\partial\rho}
{\partial r},\,\,\,\,\nabla q = h(\rho)^{1/2}\nabla \rho.
\end{equation}
Therefore (\ref{2.29}) becomes
\begin{equation}
\label{2.29c}
h(\rho)^{1/2}[\nabla^2 (q)-\frac{2}{r}\frac{\partial q}{\partial r}]=
r^2\left(S(\rho) -\Phi-\frac{H(\rho)}{r^2}\right). 
\end{equation}
If the relationship between $q(\rho)$ and $\rho$ in \eqref{2.29a} 
is invertible viz. we can express $\rho=\rho(q)$ then \eqref{2.29c}  
can be expressed in terms of $q$ only.

One possible strategy to obtain only one equation for $\rho$ is to solve
\eqref{2.29} (algebraically) for $\Phi$ and substitute the result in  
\eqref{2.4}. However, in general, this leads to a highly nonlinear equation 
for $\rho$ which has to be solved numerically.

Another way to eliminate $\Phi$ from (\ref{2.29}) is to
apply the Laplace operator to this equation and use (\ref{2.4}) to obtain one 
fourth order equation for $\rho$;
\begin{equation}
\label{2.30}
\nabla^2\left\{\frac{1}{r^2}\left[h(\rho)^{1/2}\nabla {\bf{\cdot}} ( h(\rho)^{1/2}\nabla \rho)- \frac{2h(\rho)}{r}\rho_r\right] \right\}+
4\pi G\rho= \nabla^2 S(\rho)-\nabla^2\left(\frac{H(\rho}{r^2}\right).
\end{equation}
{However, this equation is not equivalent to (\ref{2.29}). 
In fact, one can add  to the right hand side of (\ref{2.29}) a harmonic 
function $\zeta$ and still obtain (\ref{2.30}) (since $\nabla^2 \zeta=0$).
Thus (\ref{2.29}) and (\ref{2.30}) are equivalent only if can assume that 
$\zeta=0$.}
In spite of this "defect" \eqref{2.30} has the merit
of being an equation for $\rho$ only.  To use this equation to find
solutions of \eqref{2.29} and \eqref{2.4} one can implement a two stage
strategy (similar to the "predictor-corrector" in numerical analysis). 
At first one solves (\ref{2.30}) to find (general) ``initial" solution for
$\rho$, then use this solution in (\ref{2.4}) to find 
``initial" solution for $\Phi$. Substituting these solutions in 
(\ref{2.29}) one expunges those terms that are not consistent with this 
equation. As a result one obtain
expressions for $\Phi$ and $\rho$ that satisfy (\ref{2.4}) and
(\ref{2.29}). This final solution can be verified then by direct substitution 
of these expressions in these equations. However we shall {\bf not} use this
procedure in the following.

\subsection{The Interpretation of the Function $h(\rho)$}
 
The function $h(\rho)$ can be considered as a 
parameter function, which is determined by the momentum (and angular 
momentum) distribution in the fluid. From a practical point of view, the 
choice of this function determines the structure of the steady state 
density distribution. The corresponding flow field can be computed then 
a-posteriori (that is after solving for $\rho$) from the following relations; 
\begin{equation}
\label{2.26}
u=-\frac{1}{r}\sqrt{\frac{h(\rho)}{\rho}} \frac{\partial \rho}{\partial z}, \,\,\,
w=\frac{1}{r}\sqrt{\frac{h(\rho)}{\rho}} \frac{\partial \rho}{\partial r}.
\end{equation}

\subsection{The Steady State Pressure}

In order to derive (\ref{2.25}) we eliminated the pressure from equations
(\ref{2.12})-(\ref{2.13}). However, in practical astrophysical 
applications, it is important to know the equation of state of the fluid
under consideration. For this reason, we derive here an equation analogous to
(\ref{2.25}) for the steady state pressure. To this end, we divide 
(\ref{2.12})-(\ref{2.13}) by $\rho$, differentiate the first with respect
to $z$ the second with respect to $r$ and subtract. Using (\ref{2.6}) this 
leads to
\begin{equation}
\label{3.1}
-\frac{u}{r}\chi +u\frac{\partial\chi}{\partial r}+ 
w\frac{\partial\chi}{\partial z} = \frac{1}{\rho^2}J\{\rho,p\}
-J\{\rho,\frac{f\,f_{\rho}}{r^2}\}.
\end{equation}
Expressing $u,w$ and $\chi$ in terms of $\psi$ this yields
\begin{equation}
\label{3.2}
\rho^2\,J\left\{\psi,\frac{1}{r^2}\left[\nabla^2 \psi-\frac{2}{r}\frac{\partial\chi}{\partial r}\right] \right\} = J\{\rho,p \}-
\rho^2J\{\rho,\frac{f\,f_{\rho}}{r^2}\}.
\end{equation}
Eliminating $\psi$ from this equation (using (\ref{2.20})) leads to; 
\begin{equation}
\label{3.3}
J\left\{\rho,\frac{1}{r^2}\left[\rho\psi_{\rho}^2(\nabla^2\rho -\frac{2}{r}\rho_r)+\rho\psi_{\rho}\psi_{\rho\rho}(\rho_r^2+\rho_z^2)\right] \right\} 
=\frac{1}{\rho}J\{\rho,p \}-\rho J\{\rho,\frac{f\,f_{\rho}}{r^2}\}.
\end{equation}
{This equation is satisfied if there exists some function 
$P(\rho)$ such that} 
\begin{equation}
\label{3.4}
h(\rho)\left( \nabla^2 \rho -\frac{2}{r}\rho_r\right)+ \frac{1}{2}\left[h^{\prime}(\rho) - \psi_{\rho}^2 
\right] (\rho_r^2 + \rho_z^2)= r^2\left(\frac{p}{\rho} - 
\frac{\rho f\,f_{\rho}}{r^2}+P(\rho)\right).
\end{equation}
where $P(\rho)$ is some function of $\rho$. 
Subtracting this equation from (\ref{2.25}), we then have
\begin{equation}
\label{3.5}
\frac{p}{\rho}= S(\rho) - P(\rho) - \frac{1}{2r^2} \psi_{\rho}^2 
(\rho_r^2 + \rho_z^2) -\Phi -\frac{f^2}{2r^2},
\end{equation}
or equivalently
\begin{equation}
\label{3.5z}
p = \rho(S(\rho) - P(\rho)) - \frac{h(\rho)}{2r^2} 
(\rho_r^2 + \rho_z^2) -\rho(\Phi -\frac{f^2}{2r^2}).
\end{equation}

Therefore, the solution of (\ref{2.25}) and (\ref{2.4}) determines
the pressure distribution in the fluid (assuming that the functions
$P,S$ {are known}). 

Conversely, if the pressure distribution is known apriori, e.g if we assume that
the fluid is a polytropic gas where $p=A\rho^{\alpha+1}$, then (\ref{2.4})
can be used to eliminate $\Phi$ from (\ref{3.5}). 
\begin{equation}
\label{3.6}
\nabla^2(P) = \nabla^2\left[S-A\rho^{\alpha}- \frac{1}{2r^2} \psi_{\rho}^2
(\rho_r^2 + \rho_z^2) \right] - 4\pi G \rho-
\nabla^2\left(\frac{f^2}{2r^2}\right).
\end{equation}
{As in the derivation of (\ref{2.30}) from (\ref{2.25}) 
this differentiation implies that the solutions of (\ref{3.5}) and (\ref{3.6}) 
might differ by a harmonic function term.}

It follows then that for a polytropic gas eqs. (\ref{3.4}),(\ref{3.6})
form a closed system of coupled equations for $\rho$ and $P$ with 
{parameter functions $\psi_{\rho}^2$ and $f(\rho)$}. 
However, if we eliminate $P$ from these two equations, we recover (\ref{2.30}). 
\setcounter{equation}{0}
\section{On the Shape of a Rotating Star}
To simplify the algebra, as far as possible, we shall let $S(\rho)=0$
and $P(\rho)=0$. With these assumptions the explicit form of \eqref{2.4}
\eqref{2.25} and \eqref{3.5z} respectively are
\begin{equation}
\label{13.1}
\frac{\partial^2\Phi}{\partial r^2}+\frac{\partial^2\Phi}{\partial z^2}+
\frac{1}{r}\frac{\partial\Phi}{\partial r}-4\pi G\rho=0,
\end{equation}
\begin{equation}
\label{13.2}
h(\rho)\left(\nabla^2\rho -\frac{2}{r}\rho_r \right) +
\frac{h^{\prime}(\rho)}{2}(\rho_r^2+\rho_z^2) = 
-r^2\left(\Phi+\frac{H(\rho)}{r^2}\right),
\end{equation}
\begin{equation}
\label{13.3}
p=-\frac{h(\rho)}{2r^2} (\rho_r^2 + \rho_z^2) -\rho(\Phi -\frac{f^2}{2r^2}).
\end{equation}
Solving \eqref{3.2} for $\Phi$ and substituting in \eqref{3.3} we obtain
\begin{eqnarray}
\label{13.4}
&&p=-\frac{h(\rho)(\rho_r^2 +\rho_z^2)}{2r^2})+
\rho h(\rho)\frac{2r\frac{\partial^2\rho}{\partial r^2}+
2r\frac{\partial^2\rho}{\partial z^2}-
2\frac{\partial\rho}{\partial r}}{2r^3}+ \\ \notag
&&\rho\left[\frac{dh(\rho)}{d\rho}\frac{\rho_r^2 +\rho_z^2}{2r^2}+\frac{H(\rho)}{r^2}+\frac{f^2}{2r^2}\right].
\end{eqnarray}

Along the boundary $z=z(r)$ and $p=0$ (\cite{JT} p. $54$ and p. $121$). 
The explicit expression of \eqref{13.4} along the boundary is:
\begin{eqnarray}
\label{13.5}
&&-\frac{h(\rho(r,z(r))}{2r^2}\left[\frac{\partial\rho(r,z)}{\partial r}+
\frac{\partial\rho(r,z)}{\partial z}\frac{dz(r)}{dr}\right]^2+\\ \notag
&&\frac{\rho(r,z)}{2r^2}\left[\frac{dh(\rho)}{d\rho}\left(\frac{\partial\rho(r,z)}
{\partial r}+\frac{\partial\rho(r,z)}{\partial z}\frac{dz(r)}{dr}\right)^2+2H(\rho)+f^2\right]+ \\ \notag
&&\frac{\rho(r,z(r))h(\rho(r,z(r))}{r^2}\left[\frac{\partial\rho(r,z)}{\partial z}
\frac{d^2z(r)}{dr^2}+\frac{\partial^2\rho(r,z)}{\partial z^2}
\left(\frac{dz(r)}{dr}\right)^2\right]+ \\ \notag
&&\frac{\rho(r,z(r))h(\rho(r,z(r))}{r^3}\left[\left(2r\frac{\partial^2\rho(r,z)}{\partial r \partial z}-\frac{\partial\rho(r,z)}{\partial z}\right)\frac{dz(r)}{dr}+r\frac{\partial^2\rho}{\partial r^2}-\frac{\partial\rho(r,z)}{\partial r}\right]=0.
\end{eqnarray}
Where all the partial derivatives in this expression have to be evaluated 
at $z=z(r)$. This differential equation for $z(r)$  shows clearly how the 
boundary of the star depend on the asymptotic values of $\rho$ and 
its derivatives at the boundary.

\subsection{Solutions of Equation \eqref{13.5}}
Equation \eqref{13.5} requires experimental data about the asymptotic values of 
$\rho$ and it derivatives at the boundary in order to determine the actual 
shape of the star. In absence of such data we derive in the following analytic 
solutions to \eqref{13.5} under some simplifying assumptions on these
values at the boundary (which we treat as "parameters"). 
We show that with proper choice of these parameters one 
can obtain physically acceptable "star shapes".

\subsubsection{Solutions with $h = 1$, $S = 0$, $H = 0$}

When $H=0$ the solution of \eqref{2.14a} is given by $f^2=\frac{D}{\rho}$.
To obtain a closed form solution for $z(r)$ we make the assumption that all 
second order derivatives of $\rho$ along the boundary are zero except for
$\frac{\partial^2\rho}{\partial^2 r}=F$ where $F$ is a negative constant. 
Furthermore we let $\rho$ and its first order derivatives to be constant 
along the boundary viz. $\rho(r,z(r))=A$, $\frac{\partial\rho(r,z(r)}{r}=C$ 
and $\frac{\partial\rho(r,z(r)}{z}=B$ where $A,\,B\,C$ are nonzero constants. 
As a result \eqref{13.5} reduces to 
\begin{equation}
\label{13.6}
2Ar\frac{d^2z(r)}{dr^2}-Br\left(\frac{dz(r)}{dr}\right)^2-(2Cr+2A)\frac{dz(r)}{dr}
+\frac{2A F r - C^2 r - 2AC + Dr}{B}=0
\end{equation}
Assuming that $(2AF+D)$ is negative the solution of this equation 
is
\begin{equation}
\label{3.7}
z(r)=-\frac{A}{B}\ln\left\{\frac{B^2(C_1J_1(\alpha r)-C_2Y_1(\alpha r))^2}
{4\alpha^2A^2(J_0(\alpha r)Y_1(\alpha r)-Y_0(\alpha r)J_1(\alpha r))^2}\right\}-
\frac{Cr}{B}
\end{equation}
where $\alpha^2=-\frac{2AF + D}{4A^2}$, $J$, $Y$ represent Bessel functions 
of the first and second kind and $C_1$, $C_2$ are integration 
constants. The value of these constants can be determined  using 
the boundary conditions $z(0)=1$, $z(1)=0$. A reasonable shape for the radius
of star can be obtained then with $A=1$, $C=-1$, $B=-10$ and $\alpha=0.0001$. 
This shape is shown in $Fig\,\, 1$

We note that if one assigns (constant) non zero values to the second 
order derivatives of $\rho$ at the boundary one obtains a solution of 
\eqref{13.5} in terms of (lengthy) expressions with Whitaker functions.

\subsubsection{Solutions with $h = 1$, $S = 0$ and $f$ is a constant}

When $f$ is a constant in \eqref{2.14a} then $H=\frac{f^2}{2}$. Using the 
same assumptions about the second and first order derivatives of $\rho$ 
as in the previous subsection we obtain the following differential equation 
for $z(r)$,
\begin{equation}
\label{3.8}
\frac{AB}{r^2}\frac{d^2z(r)}{dr^2}-\frac{B^2}{2r^2}\left(\frac{dz(r)}{dr}\right)^2-
\frac{B}{r^3}(Cr+A)\frac{dz(r)}{dr}-\frac{C^2}{2r^2}+\frac{2A(Fr-C)}{2r^3}+
\frac{Af^2}{r^2}=0
\end{equation} 
Assuming that $(f^2+F)$ is negative the solution of this equation
is
\begin{equation}
\label{3.9}
z(r)=-\frac{A}{B}\ln\left\{\frac{B^2(C_2Y_1(\beta r)-C_1J_1(\beta r))^2}
{4A^2\beta^2[J_0(\beta r)Y_1(\beta r)-Y_0(\beta r)J_1(\beta r)]^2}\right\}-
\frac{Cr}{B}
\end{equation}
where $\beta^2=-\frac{f^2 + F}{2A}$, $J$, $Y$ represent Bessel functions
of the first and second kind and $C_1$, $C_2$ are integration
constants. The value of these constants can be determined  using
the boundary conditions $z(0)=1$, $z(1)=0$. A reasonable shape for the radius
of star can be obtained then with $A=1$, $C=-1$, $B=-9$ and $\beta=0.0001$.
The resulting shape is closely similar to the one in $Fig\,\, 1$.

\subsubsection{Solutions with $h = 4\rho^2$, $S = 0$ and $H=0$}

As in Sec. $4.1$ the assumption that $H=0$ implies that $f^2=\frac{D}{\rho}$.
Assuming that {\bf all} the second order derivatives of $\rho$ are zero
and its first order derivatives are constant at the boundary (with same 
notation as before) we obtain the following differential equation for $z(r)$,
\begin{equation}
\label{3.10}
{8A^2Br}\frac{d^2z(r)}{dr^2}+{4AB^2r}\left(\frac{dz(r)}{dr}\right)^2+
{8AB}(Cr-A)\frac{dz(r)}{dr}+\frac{Dr}{A}+{4AC(Cr-2A)}=0.
\end{equation}
The solution of this equation is
\begin{equation}
\label{3.11}
z(r)=\frac{\ln \! \left(\frac{B^{2} \left(C_{2} Y_{1}\left(\alpha  r \right) \ -C_{1} J_{1}\left(\alpha  r \right) \right)^{2}}{4 A^{2} \alpha^{2} \left(J_{0}\left(\alpha  r \right) Y_{1}\left(\alpha  r \right)-J_{1}\left(\alpha  r \right) Y_{0}\left(\alpha  r \right)\right)^{2}}\right) A -C r}{B}
\end{equation}
where $\alpha^2=\frac{D}{16A^4}$ and $C_1$, $C_2$ are integration constants.
The value of these constants can be determined  using
the boundary conditions $z(0)=1$, $z(1)=0$. A reasonable shape for the radius
of star can be obtained then with $A=1$, $C=-1$, $B=-8.8$ and $\beta=0.001$.
The resulting shape is depicted in $Fig\,\, 2$.

\subsubsection{Solutions with $h = 4\rho^2$, $S = 0$ and $f$ is a constant}

When $f$ is a constant in \eqref{2.14a} then $H=\frac{f^2}{2}$. Using the
same assumptions about the second and first order derivatives of $\rho$
as in the previous subsection we obtain the following differential equation
for $z(r)$,
\begin{equation}
\label{3.12}
{8A^2Br}\frac{d^2z(r)}{dr^2}+{4AB^2r}\left(\frac{dz(r)}{dr}\right)^2+
{8AB}(Cr-A)\frac{dz(r)}{dr}+{4AC(Cr-2A)}+2f^2r=0.
\end{equation}
The solution of this equation is formally the same as the one in \eqref{3.11}
with $\alpha^2=\frac{f^2}{8A^3}$. The plot for the radius of the star is
similar to the one in $Fig\,\,2$

As a generalization of the results in this section one might consider the
case where $H(\rho)=(n+1)\rho^n$ where $n$ is a positive integer. In this case 
$f^2=2\rho^n$. Using the same procedure described in this section one obtains
a solution for the star radius which is similar to the one presented in 
$Fig\,\,2$
\setcounter{equation}{0}
\section{Special Solutions for the Shape of a Rotating Star }
In this section we present some solutions for the shape of a axi-symmetric 
star or interstellar cloud subject to some 
assumptions regarding the functional dependence of $\rho$ on $z$.

\subsection{Solutions with $h=1,S=0,H=0$}

Eq. (\ref{2.29}) is, in general, a nonlinear equation, which (to our
best knowledge) can not be solved (in general) analytically. The only
exception is the case where $h$ is a constant under which the resulting
equation is linear. It should be remembered, however, that although
(\ref{2.29}) reduces to a linear equation when $h$ is a constant,
the original equations (\ref{2.1})-(\ref{2.4}) of the model are nonlinear
for this choice of $h$ as is evident from (\ref{2.26}). Therefore,
in principle, we are still attempting to solve to a system of nonlinear
equations.

For this choice of $h$,\,we have from (\ref{2.26}) that
\begin{equation}
\label{4.02}
(u,w) = \frac{1}{r\sqrt{\rho}}\left(-\frac{\partial \rho}{\partial z},
\frac{\partial \rho}{\partial r}\right),\,\,\,\, \rho > 0.
\end{equation}
That is with the same gradient of $\rho$, $(u,w)$
will increase as $\rho$ decreases. We conclude then that, in general,
matter in regions with low density might have higher momentum than in
regions of higher density.

In the following we consider \eqref{2.29c} and \eqref{2.4} with
$h(\rho)=1$, $H(\rho)=0$ and $S(\rho)=0$.
For this choice of $H(\rho)$ we obtain from
\eqref{2.14a} that $f^2(\rho)=\frac{C^2}{\rho}$ where $C$ is a constant.

Under this assumptions it follows from \eqref{2.29} reduces to
\begin{equation}
\label{4.03}
\nabla^2 \rho -\frac{2\rho_r}{r}=-r^2\Phi.
\end{equation}  
Thus \eqref{2.4} and \eqref{4.03} represent a system of of two coupled linear 
equations for $\rho$ and $\Phi$. To solve this system we start by applying 
separation of variables for $\rho$ and $\Phi$. Thus we set 
\begin{equation}
\label{4.03a}
\rho(r,z)=R(r)Z(z),\,\,\,\Phi(r,z)=\phi(r)\psi(z)
\end{equation}
Substituting these expressions in \eqref{2.4} and \eqref{4.03} leads to
\begin{equation}
\label{4.04}
\frac{\frac{d^2R}{dr^2}}{R} -\frac{\frac{dR}{dr}}{rR}+ 
\frac{\frac{d^2Z}{dz^2}}{Z} + r^2\frac{\phi\psi}{RZ}=0,
\end{equation}  
\begin{equation}
\label{4.05}
\frac{\frac{d^2\phi}{dr^2}}{\phi} +\frac{\frac{d\phi}{dr}}{r\phi}+ 
\frac{\frac{d^2\psi}{dz^2}}{\psi} -4G\pi\frac{RZ}{\phi\psi}=0.
\end{equation}
It is clear that in order to make further progress some functional 
relationships must be introduced between $Z(z)$ and $\psi(z)$.
In the following we consider two such relationships (and assume further that
$\rho(r,z)$ is symmetric with respect to the $x-y$ plane viz 
$rho(r,z)=rho(r,-z)$).
\subsubsection{$Z(z)=\psi(z)=D_1z+D_2$,\,\,\, $D_1\ne 0$}

Under these assumptions \eqref{4.04} and \eqref{4.05} become
\begin{equation}
\label{4.06}
\frac{\frac{d^2\phi}{dr^2}}{\phi} +\frac{\frac{d\phi}{dr}}{r\phi}-
4G\pi\frac{R}{\phi}=0.
\end{equation}
\begin{equation}
\label{4.07}
\frac{\frac{d^2R}{dr^2}}{R} -\frac{\frac{dR}{dr}}{rR}+ 
r^2\frac{\phi}{R}=0,
\end{equation}
Solving \eqref{4.07} algebraically for $\phi$  we have
\begin{equation}
\label{4.08}
\phi(r)=-\frac{r\frac{d^2R}{dr^2}- \frac{dR}{dr})}{r^3}
\end{equation}
Substituting this expression in \eqref{4.06}
we find that $R(r)$ has to satisfy the following fourth order equation.
\begin{equation}
\label{13.10}
r^3\frac{d^4R}{dr^4} - 4r^2\frac{d^3R}{dr^3} + 9r\frac{d^2R}{dr^2} - 
9\frac{dR}{dr}+4G\pi r^5R=0. 
\end{equation}
This equation has analytic solution in term of three hypergeometric functions
and a Meijer function.
\begin{eqnarray}
\label{13.11}
&&R(r)=C_1r^4F([0,3], [1, 4/3, 5/3], \alpha) + 
C_2F([0,3], [1/3, 1/3, 2/3], \alpha) + \\ \notag
&&C_3r^2F([0,3], [2/3, 2/3, 4/3], \alpha)+
C_4MeijerG([[2,0], [0,4]], [[2/3, 2/3], [1/3, 0]],\alpha)
\end{eqnarray}
where $\alpha=\frac{-\pi Gr^6}{324}$ and $C_i,\,i=1,2,3,4$ are arbitrary
constants.

To determine the shape of the resulting star we impose on the pressure
the condition $p=0$ at the boundary.(\cite{JT} p. $54$ and p. $121$).
To compute the pressure, and impose this boundary condition we use
(\ref{3.5z}) with $S(\rho)=P(\rho)=0$. Under these settings \eqref{3.5z}
becomes
\begin{equation}
\label{13.12}
z=\frac{\sqrt{D}-D_2}{D_1}
\end{equation}
where
$$
D=\frac{r(D_1^2R^2 - C^2)}{2rR\frac{d^2R}{dr^2} - r(\frac{dR}{dr})^2 - 
2R\frac{dR}{dr}}
$$
It should be observed that the choice of the parameters in this equation 
is subject to some constraints. Thus we must have 
$\rho(r,z)=R(r)(D_1z+D_2) \ge 0$ 
throughout the domain. Furthermore the value of $R(r)$ 
at the boundary of the star should be close to zero. Finally the shape of 
the star should conform to the observed physical (astronomical) data.

We simulated \eqref{13.10} (numerically) on the interval $r_0 \le r \le 1$
(with $r_0=0.001$ to avoid numerical issues at $r=0$) with the boundary 
conditions
$$
R(r_0)=1,\,\,\, \frac{dR}{dr}(r_0)= -0.003,\,\,\, R(1)=10^{-5} ,\,\,\ 
\frac{dR}{dr}(1)=-0.015.
$$ 
For the parameters we used $C=1$, $D_2=-D_1=2.543\times 10^{2} $ 
and obtained Fig $3$ for the star radius $R_s$ and $z$ as a function of $r$.

\subsubsection{$Z(z)=\psi(z)=D_3\exp{\lambda z}$}
We observe that the dependence of the system \eqref{4.04} and \eqref{4.05} 
on $z$ can be eliminated (also) by introducing the following 
ansatz
$$
\psi(z)=D_3\exp(\lambda z),\,\, Z(z)=D_4\exp(\lambda z).
$$
(For simplicity we let $D_3=D_4$ in the following).

Substituting these expressions in \eqref{4.04} and \eqref{4.05} yield
\begin{equation}
\label{4.06a}
\frac{\frac{d^2\phi}{dr^2}}{\phi} +\frac{\frac{d\phi}{dr}}{r\phi}-
4\pi G\frac{R}{\phi}+\lambda^2=0.
\end{equation}
\begin{equation}
\label{4.07a}
\frac{\frac{d^2R}{dr^2}}{R} -\frac{\frac{dR}{dr}}{rR}+ 
r^2\frac{\phi}{R}+\lambda^2=0.
\end{equation}
Solving \eqref{4.07a} (algebraically) for $\phi$ yields,
\begin{equation}
\label{4.07b}
\phi=-\frac{\lambda^2rR+\frac{d^2R}{dr^2}{R}-\frac{dR}{dr}}{r^3}.
\end{equation}
Substituting this expression in 
\eqref{4.06a} we obtain the following fourth order equation for $R$
\begin{equation}
\label{13.10a}
r^3\frac{d^4R}{dr^4} - 4r^2\frac{d^3R}{dr^3} + 
(2\lambda^2r^3+9r)\frac{d^2R}{dr^2} - (4\lambda^2r^2+9)\frac{dR}{dr}+
(4\pi Gr^5+\lambda^4r^3+4\lambda^2r)R=0. 
\end{equation}
When $G=0$ this equation has a solution in terms of Bessel functions 
of the first and second kind of order one.
Motivated by this result we look for solutions of \eqref{13.10a} in the form
\begin{equation}
\label{4.10a}
R(r)=A(r)J_1(\lambda r)+B(r)J_0(\lambda r)
\end{equation}
(where $J_0$ and $J_1$ are Bessel functions of the first kind of order
zero and one).
Substituting this expression in \eqref{13.10a} and using the fact that
$J_0$ and $J_1$ are independent we obtain a coupled system of fourth order
equations for $A(r)$ and $B(r)$. This system was solved numerically with 
$G=1$ subject to the boundary conditions
$A(0)=B(0)=1$ and $A(1)$,$B(1)$ are zero with their first and second 
order derivatives at this point.

To determine the shape of the resulting star we impose on the pressure
the condition $p=0$ at the boundary.(\cite{JT} p. $54$ and p. $121$).
To compute the pressure, and impose this boundary condition we use
(\ref{3.5z}) with $S(\rho)=P(\rho)=0$. Under these settings \eqref{3.5z}
becomes
\begin{equation}
\label{23.5}
p = -\frac{h(\rho)}{2r^2} 
(\rho_r^2 + \rho_z^2) -\rho\left(\Phi -\frac{f^2}{2r^2}\right).
\end{equation}
Using \eqref{4.03a},\eqref{4.10a}, \eqref{4.07b} and $f^2=\frac{C}{\rho}$
where $C$ is a constant in \eqref{23.5} leads to a quadratic equation for
$w=\exp(\lambda z)$ as a function of $r$. 

We solved \eqref{13.10a} and 
\eqref{4.10a} with $\lambda=-3.6$ for $R(r)$. Substituting this result in 
\eqref{23.5} and solving the quadratic equation for $w$ with 
$C=9.658\times 10^{-3}$ yields 
Figure $4$ for the star radius $R_s$ and $z$ as a function of $r$.
\subsection{Solutions with $G \ll 1$}

In this subsection we consider a {\bf diffuse gas cloud} where $G \ll 1$ 
and solve \eqref{2.29} and \eqref{2.4} with $h(\rho)=1$, $H(\rho)=0$ and 
$S(\rho)=0$. For this choice of $H(\rho)$ we obtain from 
\eqref{2.14a} that $f^2(\rho)=\frac{C^2}{\rho}$ where $C$ is a constant.

Under these settings \eqref{2.29} takes the following form
\begin{equation}
\label{4.1}
\nabla^2 \rho - \frac{2}{r}\rho_r + r^2 \Phi=0
\end{equation}
Solving this equation algebraically 
for $\Phi(r,z)$ and substituting the result in \eqref{2.4} we obtain,
\begin{equation}
\label{4.2}
\frac{1}{r^2}\rho_{rrrr}+\frac{2}{r^2}\rho_{rrzz}+\frac{1}{r^2}\rho_{zzzz}
-\frac{4}{r^3}\rho_{rrr} -\frac{4}{r^3}\rho_{rzz}+
\frac{9}{r^4}\rho_{rr}+\frac{4}{r^4}\rho_{zz}-\frac{9}{r^5}\rho_r+4\pi G\rho=0.
\end{equation}

To this equation we can apply separation of variables viz $\rho(r,z)=R(r)Z(z)$
by introducing the ansatz that 
\begin{equation}
\label{4.3}
\frac{d^2Z}{dz^2}=\lambda^2 Z(z).
\end{equation}
The resulting equation for $R(r)$ is
\begin{equation}
\label{4.4}
r^3\frac{d^4R}{dr^4}+4r^2\frac{d^3R}{dr^3}-r(2\lambda^2r^2+9)\frac{d^2R}{dr^2}+
9\frac{dR}{dr}-(4\pi Gr^5+\lambda^4r^3+4\lambda^2r)R=0
\end{equation} 

To solve this equation we use first order perturbation in $G$ viz.
we let
\begin{equation}
\label{4.5}
R(r)=R_0(r)+GR_1(r)+O(G^2).
\end{equation}

The equations for $R_0$ and $R_1$ respectively are
\begin{equation}
\label{4.6}
r^3\frac{d^4R_0}{dr^4} - 4r^2\frac{d^3R_0}{dr^3} + (2r^3\lambda^2+9r)\frac{d^2R_0}{dr^2} -
(4r^2\lambda^2+9)\frac{dR_0}{dr}+(r^3\lambda^4+ 4r\lambda^2)R_0=0
\end{equation}
\begin{equation}
\label{4.7}
r^3\frac{d^4R_1}{dr^4} - 4r^2\frac{d^3R_1}{dr^3} + (2r^3\lambda^2+9r)\frac{d^2R_1}{dr^2} -
(4r^2\lambda^2+9)\frac{dR_1}{dr}+(r^3\lambda^4+ 4r\lambda^2)R_0+4r^5\pi R_0=0
\end{equation}

We consider two cases.

\begin{enumerate}

\item $\lambda=0$,\, $Z(z)=Az+B$,\,\,\ $A\ne 0$)
The solution of \eqref{4.6} is
\begin{equation}
\label{4.8}
R_0=D_1 + D_2r^2 + D_3r^4 + D_4r^4\ln(r)
\end{equation}
and $R_1$ turns out to be
\begin{eqnarray}
\label{4.9}
R_1&=&-\frac{\pi(120D_4\ln(r) + 12D_3 - 67D_4}{86400}
-\frac{\pi D_2r^8}{192} - \frac{\pi D_1r^6}{24}+ \\ \notag
&&(D_7\ln(r)-\frac{D_7}{4}+D_6)\frac{r^4}{4}+\frac{D_5r^2}{2} + D_8
\end{eqnarray}
where $D_i$, $i=1,\ldots,8$ are constants.
We observe, however, that the values of the constants in these expressions
are constrained by the requirement that $\rho$ must satisfy $\rho(r,z)\ge 0$

To determine the shape of the corresponding gas cloud we substitute 
these expressions in \eqref{3.5z} to obtain (up to $O(G^2)$) a linear  
equation for $w=Az+B$ (with coefficients dependent on $r$) which determine 
$z$ as a function of $r$. 

As a simplified case we consider a solution for $\rho$ (and $\phi$ from 
\eqref{4.1}) with
$$
D_1=1,\,\,\,D_2=1,\,\,\,A=3,\,\,\,B=0,\,\,\,C=3,\,\,\, G=0.1
$$
and $D_i=0$ $i=3,\ldots,8$. We Observe that these values for the parameters 
insure that $\rho(r,z)\ge 0$.

With this set of parameters the equation for $w$ is
$$
\left[\frac{D_2}{2} +\frac{D_1}{2r^2} - G\pi r^4\left(\frac {D_2r^2}{384} + 
\frac{D_1}{48}\right)\right]f^2 - 
\left[2D_2^2 + G\pi r^2\left(D_1^2+\frac{D_2^2r^4}{6} +\frac{3D_1D_2r^2}{4}\right)\right]w=0
$$
The cloud radius $R_c$ and $z$ as a function of $r$ are depicted in Fig $5$. 

\item $\lambda\ne 0$ ($Z(z)=E\exp(\lambda z)$)
The general solution for $R_0$ is 
\begin{equation}
\label{4.10}
R_0 = (C_1r^3+C_2r)J_1(\lambda r)+(C_3r^3+C_4r)Y_1(\lambda r)
\end{equation}
where $J_1$ and $Y_1$ are Bessel functions of the first and second kind 
of order $1$ and $C_i$ $i=1..4$ are constants. 

The solution for $R_1$ consists of integrals of Bessel functions 
which can be evaluated only numerically. Assuming the cloud is 
{\bf highly diffuse}
one can neglect $R_1$ to obtain a zero order approximation for its shape. 
Substituting $C_3=C_4=0$ in \eqref{4.10} lead to the following equation for $z$ 
\begin{eqnarray}
\label{4.11}
&&\frac{C^2}{r^2}-E^2e^{2\lambda z}\left[r^4C_1^2\lambda^2 + (2C_1C_3\lambda^2 + 4C_1^2)r^2 + C_3^2\lambda^2)\right]J_1(\lambda r)^2+ \\ \notag
&&E^2e^{2\lambda z}\left[4r\lambda C_1(C_1r^2 + C_3)J_0(\lambda r)J_1(\lambda r) - \lambda^2(C_1r^2 + C_3)^2J_0(\lambda r)^2\right]=0
\end{eqnarray}
using the following values for the parameter in this equation
$$
\lambda=2,\,\,\,C_1=-C_3=0.001,\,\,\,C=1.265,\,\,\,E=10,\,\,\, G=0.1
$$
and $0 \leq r \leq 6.6$  we obtain Fig $6$ for the cloud radius $R_c$
and $z$ as a function of $r$.
\end{enumerate}

\section{Summary and Conclusions}

In this paper we considered the steady state Euler-Poisson equations 
with rotations under the additional assumption of density stratification. 
The governing equations of this model consist of six nonlinear partial 
differential equations. We showed, however, that this set of equations can 
be reduced under the assumption of axial-symmetry to two. 

We derived also a separate equations for the pressure in the star with 
special consideration for those stars composed of a polytropic fluid. 

Using this reduction and the boundary condition $p=0$ we derived a differential
equation for the shape of the boundary that depends on the asymptotic
values of $\rho$ and it derivatives at the boundary. In absence of data on 
these quantities we solved this equation in closed form to obtain 
"physically reasonable star shape" by proper choice for these quantities.

Additional special solutions for the shape of a rotating star or a cloud
were obtained by making proper assumptions about the density distribution 
within the star.

This paper does not provide a general solution to the original
classical star model described by the compressible Euler-Poisson equations. 
However, it does provide insights and analytic solutions for a subclass of 
stars described by this model.

\vspace{03in}
\centerline{\bf Declaration of conflicting interest}

The authors declare that they have no known conflicting financial
interests or personal relationships that could have appeared to influence 
the work reported in this paper

\newpage
\begin{figure}[ht!]
\includegraphics[scale=1,height=160mm,angle=0,width=180mm]{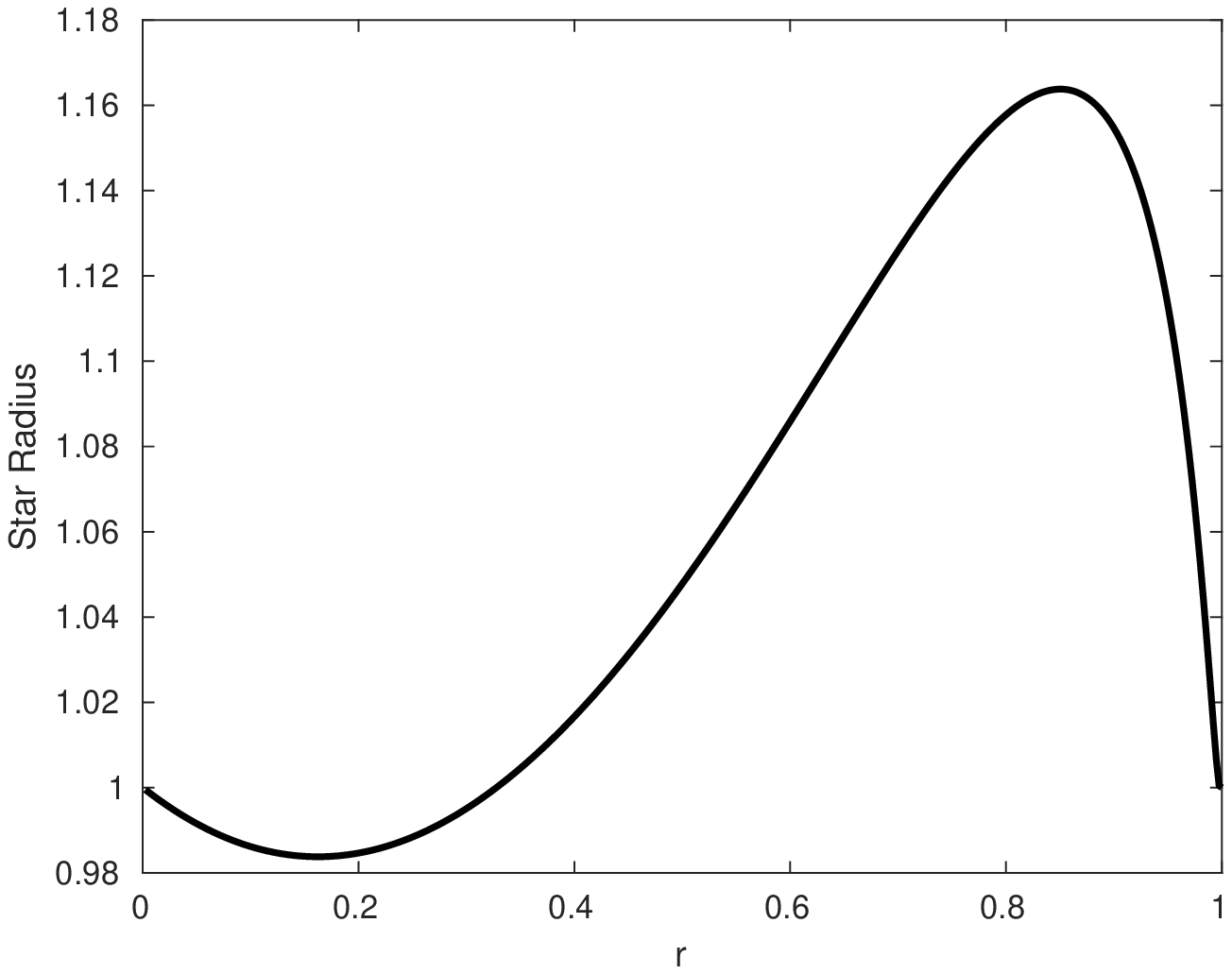}
\label{Figure 1}
\caption{Star Radius as a function of $r$ using the parameters in $Sec\,4.1.1$}
\end{figure}
\newpage
\begin{figure}[ht!]
\includegraphics[scale=1,height=160mm,angle=0,width=180mm]{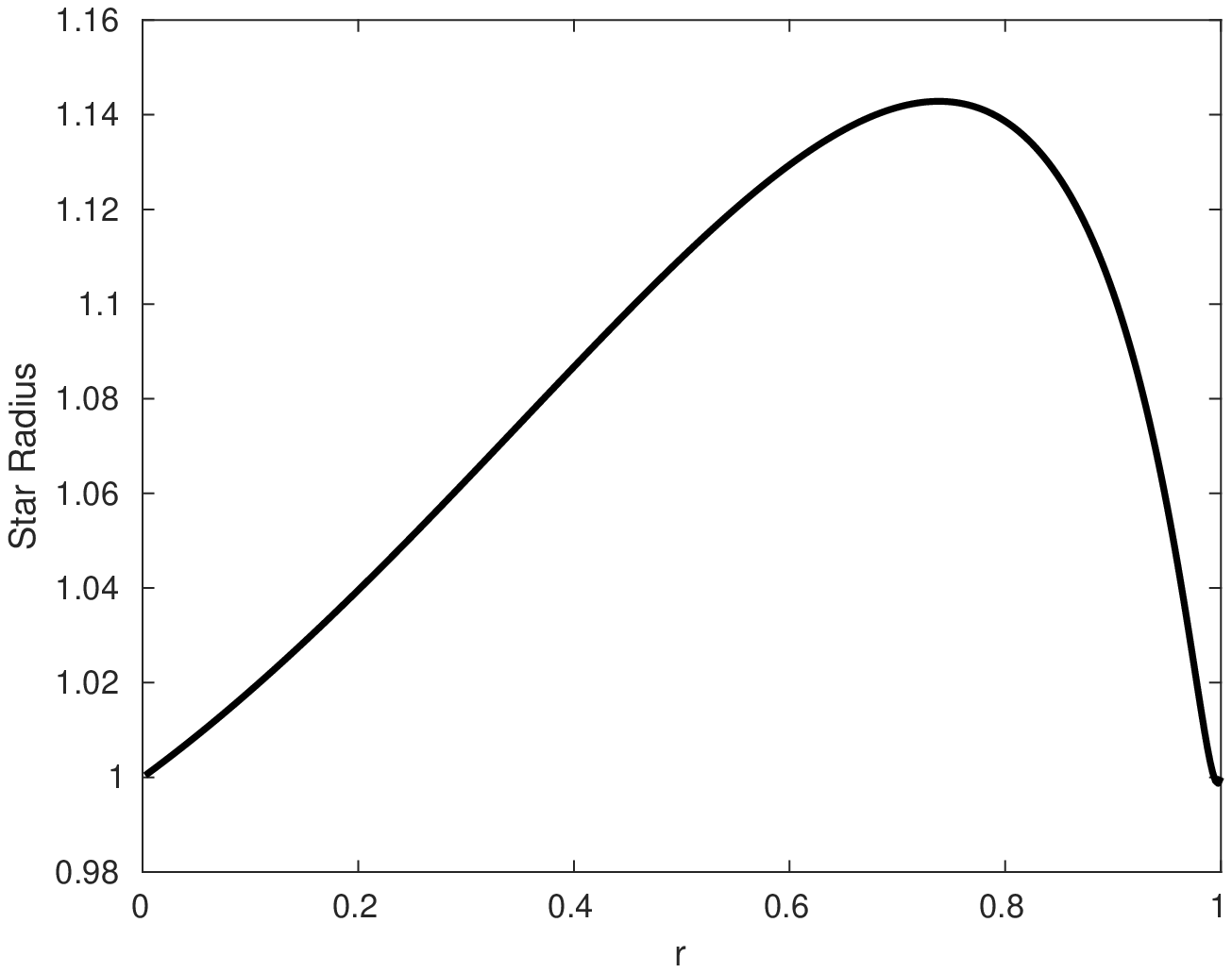}
\label{Figure 2}
\caption{Star Radius as a function of $r$ using the parameters in $Sec\,4.1.3$}
\end{figure}
\newpage
\begin{figure}[ht!]
\includegraphics[scale=1,height=160mm,angle=0,width=180mm]{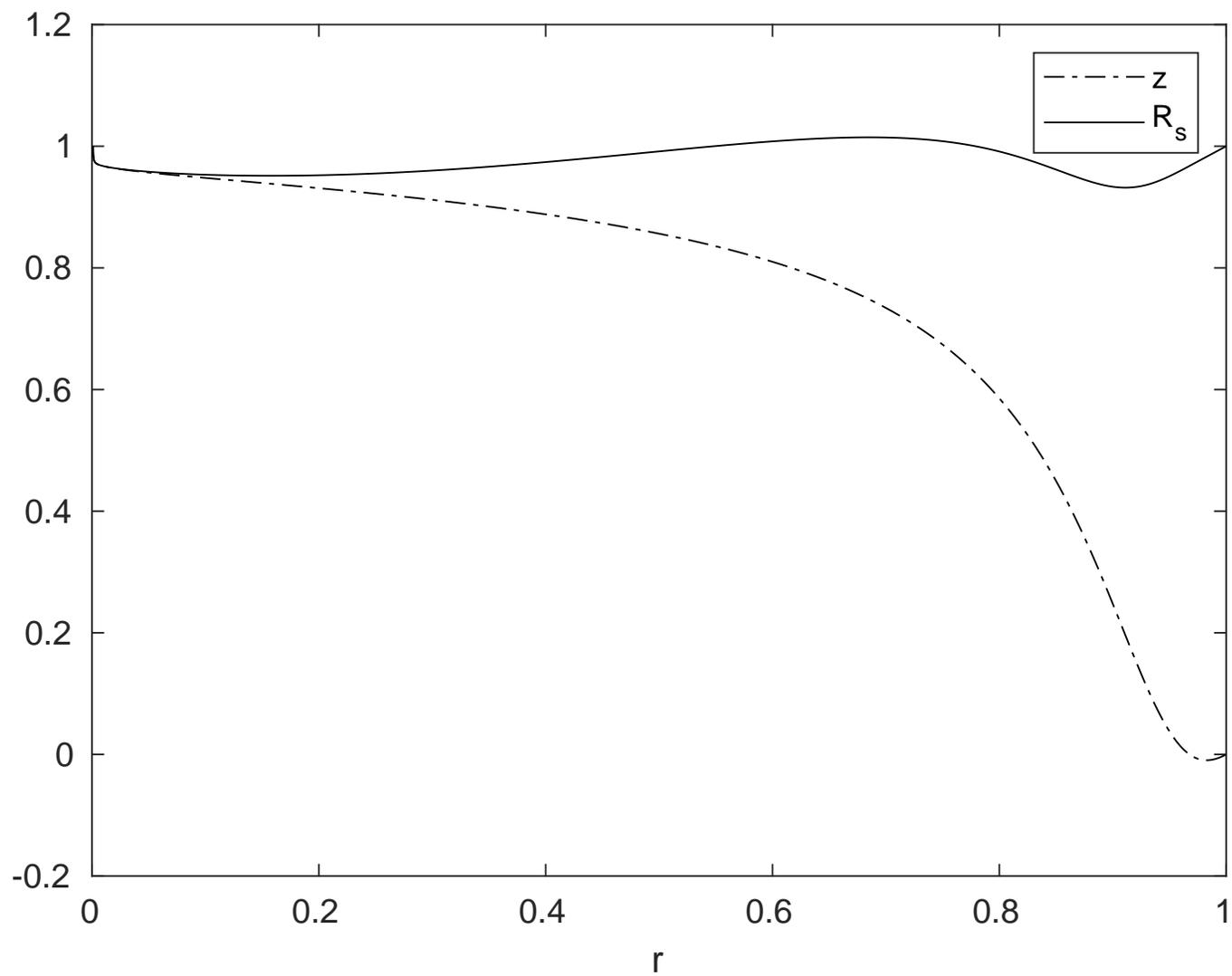}
\label{Figure 3}
\caption{Radius $R_s$ and  $z$ of a star as a function of $r$}
\end{figure}
\newpage
\begin{figure}[ht!]
\includegraphics[scale=1,height=160mm,angle=0,width=180mm]{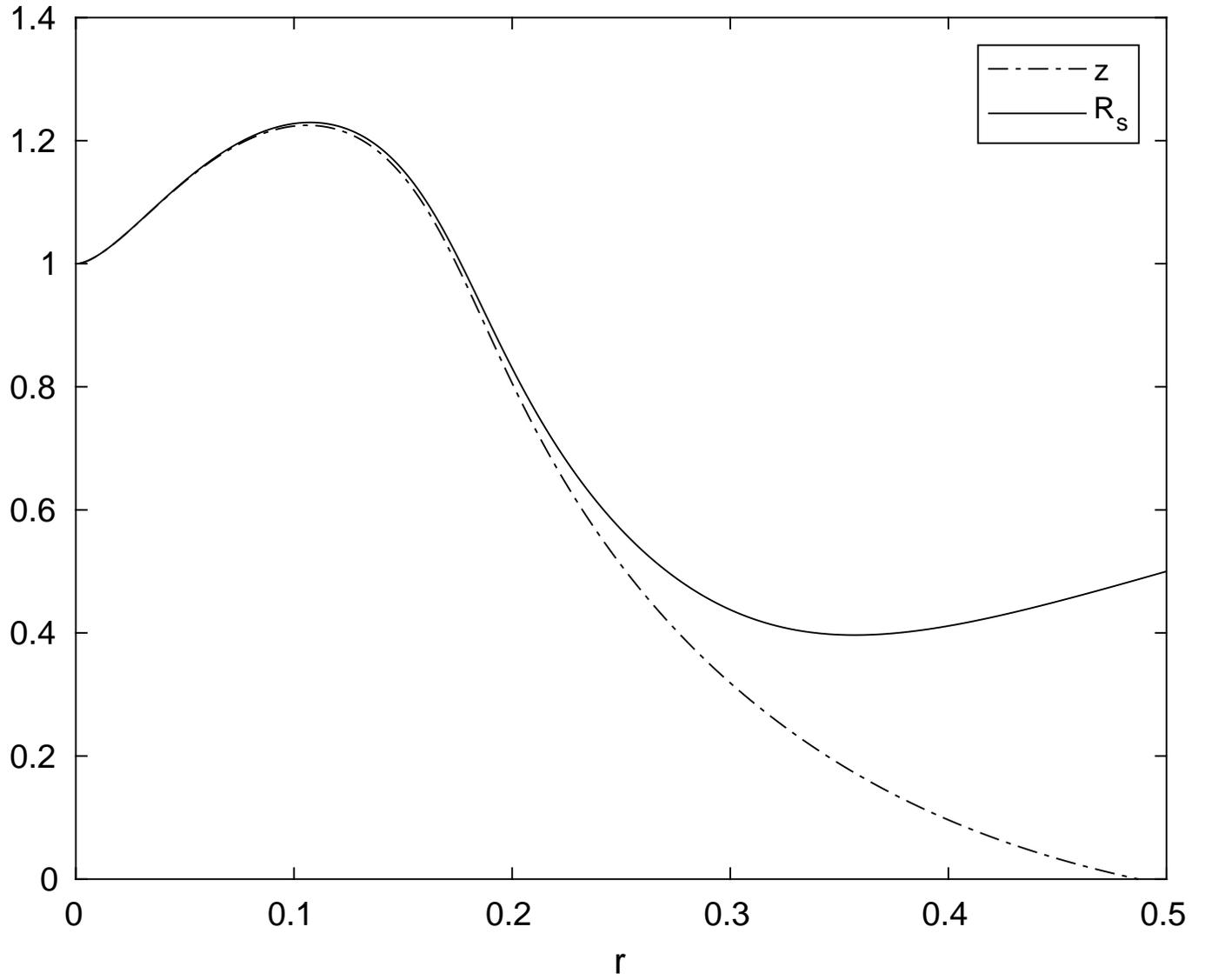}
\label{Figure 4}
\caption{Radius $R_s$ and $z$ of a star as a function of $r$, 
$\lambda=-3.6,\, C=9.658\times 10^{-3}$}
\end{figure}
\newpage
\begin{figure}[ht!]
\includegraphics[scale=1,height=160mm,angle=0,width=180mm]{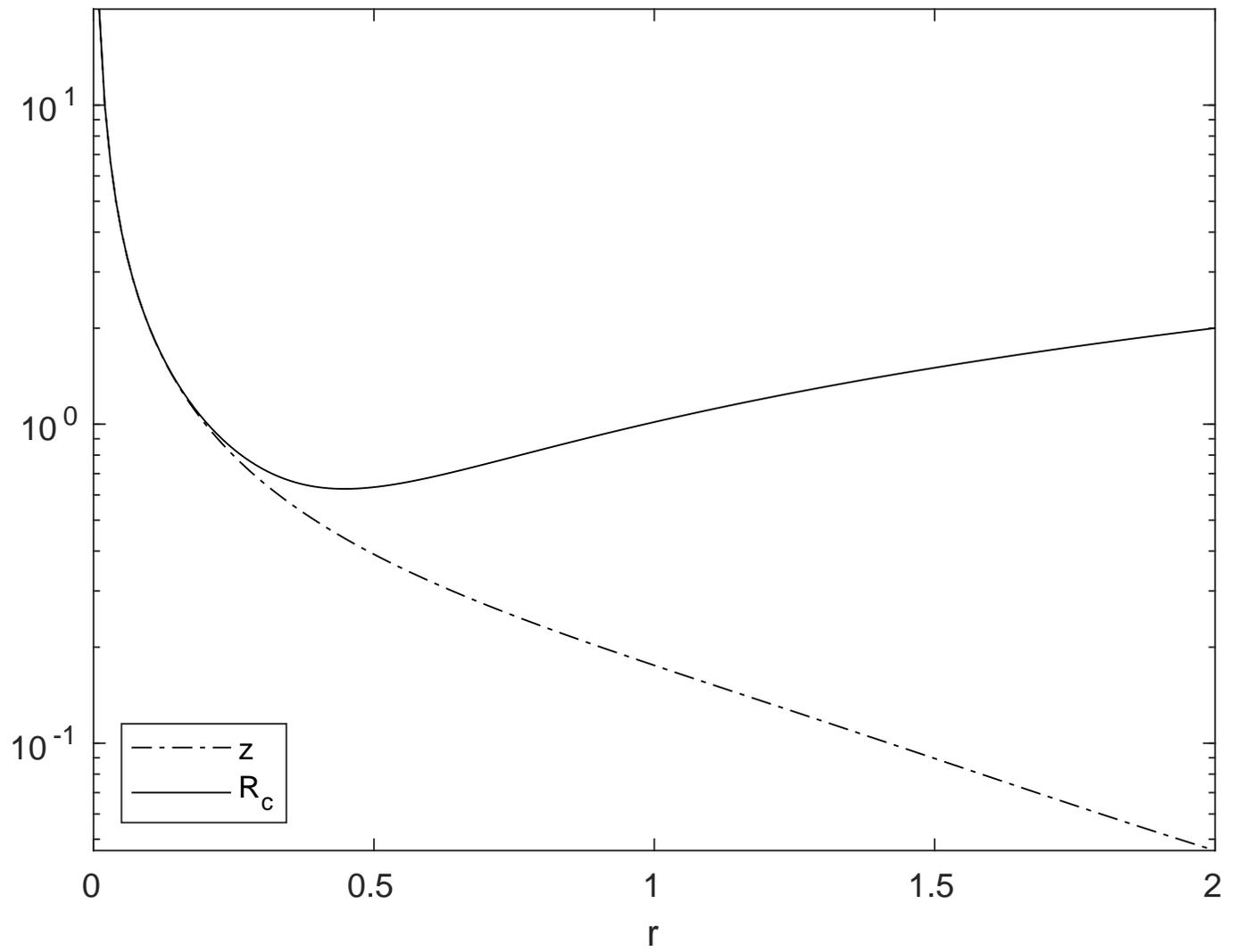}
\label{Figure 5}
\caption{$R_c$ and $z$ as a function of $r$ for a cloud with $Z(z)=Az+B$}
\end{figure}
\newpage
\begin{figure}[ht!]
\includegraphics[scale=1,height=160mm,angle=0,width=180mm]{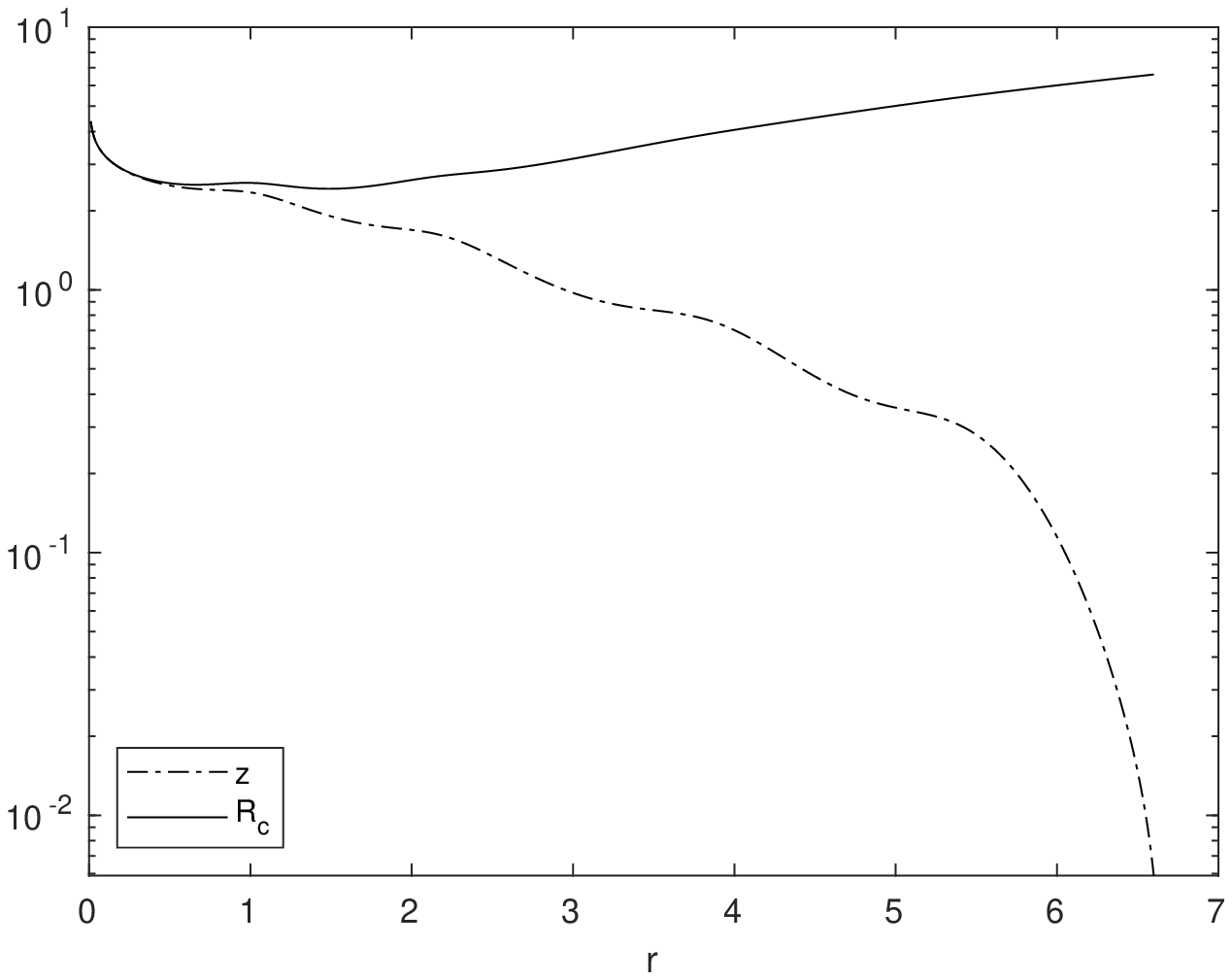}
\label{Figure 6}
\caption{$R_c$ and $z$ as a function of $r$ for a cloud with $\lambda=2$}
\end{figure}


\newpage
\begin{thebibliography}{7}
\bibitem{VSB} Beskin V. S.: Axisymmetric stationary flows in compact 
astrophysical objects, Physics-Uspekhi 40 pp. 659-688, 1997.
 
\bibitem{SC1} Chandrasekhar, S.: An Introduction to the Study of Stellar 
Structures, University of Chicago Press, Chicago, 1938.

\bibitem{SC2} Chandrasekhar, S: Ellipsoidal figures of equilibrium-
an historical account, Comm. Pure Appl. Math, 20 251-265, 1967.

\bibitem{AT} Friedman, A., Turkington, B.: The oblateness of an axisymmetric 
rotating fluid, Indiana Univ. Math fluid J., 29 pp.777-792, 1980.

\bibitem{MH1} Humi, M.: Steady States of self gravitating incompressible
fluid. J. Math. Phys. 47, 093101 (10 pages), 2006.

\bibitem{MH2} Humi, M.: Steady States of Self Gravitating Incompressible 
Fluid with Axial Symmetry, Int. J. Mod. Phys. A 24, No. 23 pp. 4287-4303 (2009)

\bibitem{MH3} Humi, M.: A Model for Pattern Formation Under Gravity, 
Applied Mathematical Modeling 40, pp. 41-49, 2016

\bibitem{MH4} Humi, M :Patterns Formation in a Self-Gravitating Isentropic Gas,
Earth Moon Planets Volume 121, Issue 1 (2018) , Pages 1-12,
https://doi.org/10.1007/s11038-017-9512-y.

\bibitem{MH5} Humi, M. and Roumas, J :Structure of polytropic stars in 
General Relativity, Astrophysics and Space Science 364(7)117, 2019
DOI: 10.1007/s10509-019-3608-y

\bibitem{MH6} Humi, M.: On the Evolution of a Primordial Interstellar Gas 
Cloud, Journal of Mathematical Physics 61, 093504 (2020);
https://doi.org/10.1063/1.5144917

\bibitem{JJM} Janga,J and Makino, T: On rotating axisymmetric solutions of 
the Euler–Poisson equations, Journal of Differential Equations
266, pp. 3942-3972, 2019

\bibitem{KNM} M. Kiguchi M., S. Narita S., Miyama S. M., Hayashi C.: 
The Equilibria of Rotating Isothermal Clouds, Astrophysical J., 
317 pp.830-845, 1987.

\bibitem{AK} Kovetz, A.: Slowly rotating polytropes. Astrophys. J. 154, 
pp. 999-1003, 1968.

\bibitem{KN} Kunzle H.P., Nester J.M.: Hamiltonian formulation of 
gravitating perfect fluids and the Newtonian limit, J. Math. Phys 25, 
pp. 1009-1018, 1984.

\bibitem{LO} Letelier, P.S., Oliveira, S.R.:  Exact self-gravitating 
disks and rings: A solitonic approach ,J. Math. Phys. 28 pp.165-170, 1987.

\bibitem{YYL} Li, Y.Y.: On uniformly rotating stars, Arch. Ration. Mech. Anal,
115 pp.367-393, 1991.

\bibitem{LS} Luo T and Smoller J.:
Existence and Non-linear Stability of Rotating Star Solutions of the 
Compressible Euler-Poisson Equations, Arch. Ration. Mech. Anal. 
191, pp.447-496, 2009.

\bibitem{MHT} Matsumoto, T. and Hanawa T.: Bar and Disk Formation 
in Gravitationally Collapsing Clouds. Astrophys. J., 521(2), pp.659-670,1999.

\bibitem{EAM} Milne, E.A : The equilibrium of a rotating star. Mon. Not. R. 
Astron. Soc. 83, pp.118-147, 1923.

\bibitem{OVM} Ortega V. G., Volkov E. and Monte-Lima I.:  Axisymmetric 
instabilities in gravitating disks with mass spectrum, Astronomy and 
Astrophysics 366, pp.276-280, 2001.

\bibitem{VPS} Paschalidis V., Stergioulas N.:
Relativistic Rotating Stars, Living Reviews in Relativity, 20, 
Article number: 7 (2017) DOI  - 10.1007/s41114-017-0008-x

\bibitem{AP1} Prentice, A. J. R., Origin of the Solar system,
Earth Moon and Planets, 19, pp. 341-398. (1978)

\bibitem{AP2} Prentice, A.J.R.; Dyt, C.P. . "A numerical simulation 
of supersonic turbulent convection relating to the formation of the 
Solar system", Monthly Notices of the Royal Astronomical 
Society. 341 (2): 644–656 (2003).

\bibitem{IWR} Roxburgh I. W., Non-Uniformly Rotating, Self-Gravitating,
Compressible Masses with Internal Meridian Circulation,
Astrophysics and Space Science, 27,pp.425-435, 1974.

\bibitem{JT} Tassoul Jean-Louis: Theory of Rotating Stars, 
Princeton U press, Princeton, NJ. 1978,

\bibitem{CSY} Yih C-S: Stratified flows. Academic Press, New York,
NY, 1980.

\end{thebibliography}
\end{document}